\documentclass{article}
\usepackage{stmaryrd}
\usepackage{mathrsfs}
\usepackage[centertags]{amsmath}
\usepackage{amsfonts, dsfont}
\usepackage{amssymb}
\usepackage{amsthm}

\input xy
\xyoption{all}

\theoremstyle{plain}
\newtheorem{thm}{Theorem}[section]
\newtheorem{cor}[thm]{Corollary}
\newtheorem{lem}[thm]{Lemma}
\newtheorem{prop}[thm]{Proposition}

\theoremstyle{definition}
\newtheorem{defn}[thm]{Definition}
\newtheorem{remark}[thm]{Remark}
\newtheorem*{ack}{Acknowledgments}

\newcommand{\bd}{\begin{defn}}
\newcommand{\ed}{\end{defn}}
\newcommand{\bl}{\begin{lem}}
\newcommand{\el}{\end{lem}}
\newcommand{\bp}{\begin{prop}}
\newcommand{\ep}{\end{prop}}
\newcommand{\bt}{\begin{thm}}
\newcommand{\et}{\end{thm}}
\newcommand{\bc}{\begin{cor}}
\newcommand{\ec}{\end{cor}}
\newcommand{\br}{\begin{remark}}
\newcommand{\er}{\end{remark}}
\newcommand{\bdi}{\begin{diagram}}
\newcommand{\edi}{\end{diagram}}
\newcommand{\beq}{\begin{equation}}
\newcommand{\eeq}{\end{equation}}
\newcommand{\ba}{\begin{array}}
\newcommand{\ea}{\end{array}}
\newcommand{\bpf}{\begin{proof}}
\newcommand{\epf}{\end{proof}}

\newcommand{\R}{\mathds{R}}
\newcommand{\Z}{\mathds{Z}}
\newcommand{\Q}{\mathds{Q}}
\newcommand{\Zp}{\mathds{Z}_{p}}
\newcommand{\Qp}{\mathds{Q}_{p}}
\newcommand{\al}{\alpha}
\newcommand{\Ga}{\Gamma}

\newcommand{\La}{\Lambda}
\newcommand{\la}{\lambda}

\newcommand{\Op}{\mathcal{O}}
\newcommand{\ord}{\mathrm{ord}}
\newcommand{\m}{\mathfrak{m}}
\newcommand{\M}{\mathfrak{M}}

\DeclareMathOperator{\Sel}{Sel} \DeclareMathOperator{\Gal}{Gal}
\DeclareMathOperator{\Hom}{Hom} \DeclareMathOperator{\rank}{rank}

\newcommand{\ot}{\otimes}
\newcommand{\ilim}{\displaystyle \mathop{\varinjlim}\limits}
\newcommand{\plim}{\displaystyle \mathop{\varprojlim}\limits}

\newcommand{\coker}{\mathrm{coker}\,}

\newcommand{\cyc}{\mathrm{cyc}}
\newcommand{\cts}{\mathrm{cts}}

\newcommand{\lra}{\longrightarrow}
\newcommand{\tha}{\twoheadrightarrow}

\newcommand{\ps}[1]{\llbracket #1 \rrbracket}

\begin{document}

\title{Comparing the Selmer group of a $p$-adic
representation and the Selmer group of the Tate dual of the
representation}
\author{Meng Fai Lim\footnote{Department of Mathematics, University of Toronto, 40 St. George St.,
Toronto, Ontario, Canada M5S 2E4}}
\date{}
\maketitle

\begin{abstract} \footnotesize
\noindent
 For a given ``ordinary" $p$-adic representation, we compare
its Selmer group with the Selmer group of its Tate dual over an
admissible $p$-adic Lie extension. Namely, we show that the
generalized Iwasawa $\mu$-invariants associated to the Pontryagin
dual of the two said Selmer groups are the same.



\end{abstract}

\section{Introduction}
The main conjecture of Iwasawa theory is a conjecture on the
relation between a Selmer group, which is a module over an (not
necessarily commutative) Iwasawa algebra, and a conjectural $p$-adic
$L$-function. This $p$-adic $L$-function is expected to satisfy a
conjectural functional equation in a certain sense. In \cite{FK},
Fukaya and Kato were able to construct such an $p$-adic $L$-function
and established a corresponding functional equation for the said
$p$-adic $L$-function assuming the (local and global) noncommutative
Tamagawa number conjecture. In view of the main conjecture and this
functional equation, one would expect to have certain algebraic
relationship between the Selmer group attached to a Galois
representation and the Selmer group attached to the Tate twist of
the dual of the Galois representation which can be thought as an
algebraic manifestation of the functional equation. It is precisely
a component of this algebraic relationship that this paper aims to
investigate.  For a module over an Iwasawa algebra, Howson and
Venjakob independently developed the notion of a generalized
$\mu$-variant which extends the classical $\mu$-invariant. The main
conjecture predicts that this invariant attached to a suitable
Selmer group should contribute a certain power of $p$ occurring in
the leading term of the conjectural $p$-adic $L$-function. In this
article, we will show that the Selmer group attached to a Galois
representation and the Selmer group attached to the Tate twist of
the dual representation have the same generalized $\mu$-variant.

In the case of a cyclotomic $\Zp$-extension, this study on the
$\mu$-invariants has been undertaken in \cite{G89, Mat}. (Actually,
in \cite{G89}, Greenberg also established the full ``algebraic"
functional equation of the Selmer groups, which we will not treat in
this article.) There he obtained the equality of the
$\mu$-invariants by combining a local-global Euler characteristic
argument with the asymptotic formulas arising from the structure
theory of $\Zp\ps{\Ga}$-modules, where $\Ga$ is the Galois group of
the cyclotomic $\Zp$-extension. When the $p$-adic Lie extension is a
multiple $\Zp$-extension, we adopt the approach of Greenberg, and in
our situation, we will need to combine the local-global Euler
characteristic argument with the asymptotic formula of Cucuo and
Monsky \cite{CM, Mon} which allow us to deal with
$\Zp\ps{G}$-modules for $G \cong \Zp^r$. We also mention that since
there are infinite decomposition of primes in a $\Zp^r$-extension
when $r\geq 2$, we will need a slightly more careful argument.

When the $p$-adic Lie extension is not commutative, the above
approach breaks down, as one does not have an asymptotic formula as
in the commutative situation, However, if we assume further that our
Selmer group satisfies a certain torsion property, then we can
establish the equality of the $\mu$-invariants via a different (and
a rather indirect) approach. Namely, under the assumption of the
above said torsion property (with some other assumptions, see
Theorem \ref{MHG} for details), we show that the $\mu$-invariant of
the Selmer group over the $p$-adic extension coincides with
$\mu$-invariant of the said Selmer group over the cyclotomic
$\Zp$-extension. This in turn allows us to deduce the equality of
the $\mu$-invariants in the noncommutative setting from the
cyclotomic $\Zp$-extension case.

We like to emphasize that although this paper is highly motivated by
the main conjecture of Iwasawa and the functional equation of the
conjectural $p$-adic $L$-function, we do not assume these
conjectures (other than the torsion property on the Selmer groups in
the noncommutative $p$-adic Lie extension situation) in all our
argument.

We now give a brief description of the layout of the paper. In
Section \ref{Algebraic Preliminaries}, we introduce the generalized
$\mu$-invariant of Howson and Venjakob, and record certain estimates
on the order of certain cohomology groups. In Section \ref{a ratio
formula}, we recalled a ratio formula of Greenberg for the strict
Selmer group of the dual of a finite module and the strict Selmer
group of the Tate twist of the Pontryagin dual of the finite module.
We will combine this ratio formula with the asymptotic formulas of
Cucuo-Monsky to equate the $\mu$-invariant of the strict Selmer
group of a torsion Galois module and the $\mu$-invariant of the
strict Selmer group of the Tate twist of its dual over a
$\Zp^r$-extension in Section \ref{compare}. In Section
\ref{compare2}, we carry out the study over a general noncommutative
$p$-adic Lie extension. In Section \ref{Artin twist of Selmer
groups}, we compare the strict Selmer groups of the Artin twists of
the representations. In Section \ref{Greenberg Selmer groups}, we
compare the strict Selmer groups with another Selmer groups of
Greenberg and the Selmer complexes. Via these comparison, we see
that the conclusion in the Greenberg strict Selmer groups can be
carried over to these groups and complexes. In Section
\ref{examples}, we finally discuss some examples of Galois
representations where our results can be applied.

\begin{ack}
     This work was written up when the author is a Postdoctoral fellow at the GANITA Lab
    at the University of Toronto. He would like to acknowledge the
    hospitality and conducive working conditions provided by the GANITA
    Lab and the University of Toronto.
        \end{ack}

\section{Algebraic Preliminaries} \label{Algebraic Preliminaries}

In this section, we recall some algebraic preliminaries that will be
required in the later part of the paper. Fix a prime $p$. Let $\Op$
be the ring of integers of some finite extension $K$ of $\Qp$. We
fix a local parameter $\pi$ for $\Op$ and denote the residue field
of $\Op$ by $k$. Let $G$ be a compact pro-$p$ $p$-adic Lie group
without $p$-torsion. It is well known that $\Op\llbracket
G\rrbracket$ is an Auslander regular ring (cf. \cite[Theorems
3.26]{V02}). Furthermore, the ring $\Op\ps{G}$ has no zero divisors
(cf.\ \cite{Neu}), and therefore, admits a skew field $Q(G)$ which
is flat over $\Op\ps{G}$ (see \cite[Chapters 6 and 10]{GW} or
\cite[Chapter 4, \S 9 and \S 10]{Lam}). If $M$ is a finitely
generated $\Op\ps{G}$-module, we define the $\Op\ps{G}$-rank of $M$
to be
$$ \rank_{\Op\ps{G}}M  = \dim_{Q(G)} Q(G)\ot_{\Op\ps{G}}M. $$
For a general compact $p$-adic Lie group $G$, we follow \cite{CH}
and extend the definition of the $\Op\ps{G}$-rank by the formula
\[ \rank_{\Op\ps{G}} N=
 \displaystyle\frac{\rank_{\Op\ps{G_0}}N}{[G:G_0]}, \]
 where $G_0$ is an open normal uniform pro-$p$ subgroup of $G$.

\bl \label{rank lemma}
 The above definition for $\Op\ps{G}$-rank is independent of the choice of
 $G_0$.
\el

\bpf Let $G_1$ be another open normal uniform pro-$p$ subgroup of
$G$. Since $G_0\cap G_1$ is also an open normal uniform pro-$p$
subgroup of $G$, one is therefore reduced to proving the equality
\[ \rank_{\Op\ps{G_1}}M = [G_0: G_1]\rank_{\Op\ps{G_0}}M\] whenever
$G_1\subseteq G_0$. Fix a finite free $\Op\ps{G_0}$-resolution
 \[ 0\lra \Op\ps{G_0}^{n_d}
\lra \cdots \lra \Op\ps{G_0}^{n_0}\lra M \lra 0\] of $M$. Then the
groups $H_i(G_0,M)$ can be computed by the homology of the complex
 \[\Op^{n_d}\lra \cdots \lra \Op^{n_0},\]
 and consequently, we obtain
 \[ \sum_{i=0}^d (-1)^i\rank_{\Op}H_i(G_0,M)=
 \sum_{i= 0}^d (-1)^i n_i.\]
 On the other hand, the above $\Op\ps{G_0}$-free resolution is
also a $\Op\ps{G_1}$-free resolution for $M$. Therefore, the groups
$H_i(G_1,M)$ can be computed by the homology of the complex
 \[  \Op^{[G_0:G_1]n_d}\lra \cdots \lra \Op^{[G_0:G_1]n_0} \]
  which gives
 \[ \ba{rl}
  \displaystyle \sum_{i=0}^d (-1)^i\rank_{R}H_i(G_1,M)\!\!\! &=  \displaystyle[G_0:G_1]\sum_{i= 0}^d (-1)^i
 n_i \\
 &= \displaystyle [G_0:G_1]\sum_{i=0}^d (-1)^i\rank_{\Op}H_i(G_0,M)
 \ea \]
 Applying a formula of Howson (cf. \cite[Theorem 1.1]{Ho}), we have
\[ \rank_{\Op\ps{G_1}}M = [G_0: G_1]\rank_{\Op\ps{G_0}}M\]
 as required.
\epf

Note that the $\Op\ps{G}$-rank needs not be an integer in general.
We will say that a $\Op\ps{G}$-module $M$ is \textit{torsion} if
$\rank_{\Op\ps{G}} M = 0$.

Now suppose that $N$ is a $k\ps{G}$-module. We then define its
$k\ps{G}$-rank by
 $$ \rank_{k\ps{G}} N=
 \displaystyle\frac{\rank_{k\ps{G_0}}N}{[G:G_0]}, $$
where $G_0$ is an open normal uniform pro-$p$ subgroup of $G$. By a
similar argument as above, one can show that this definition is
independent of the choice of $G_0$ (see also \cite[Proposition
1.6]{Ho}). Similarly, we will say that that the module $N$ is a
\textit{torsion} $k\ps{G}$-module if $\rank_{k\ps{G}}N = 0$.

For a given finitely generated $\Op\ps{G}$-module $M$, we denote
$M(\pi)$ to be the $\Op\ps{G}$-submodule of $M$ which consists of
elements of $M$ that are annihilated by some power of $\pi$. Since
the ring $\Op\ps{G}$ is Noetherian, the module $M(\pi)$ is finitely
generated over $\Op\ps{G}$. Therefore, one can find an integer
$r\geq 0$ such that $\pi^r$ annihilates $M(\pi)$. Following
\cite[Formula (33)]{Ho}, we define
  \[\mu_{\Op\ps{G}}(M) = \sum_{i\geq 0}\rank_{k\ps{G}}\big(\pi^i
   M(\pi)/\pi^{i+1}\big). \]
(For another alternative, but equivalent, definition, see
\cite[Definition 3.32]{V02}.) By the above discussion and our
definition of $k\ps{G}$-rank, the sum on the right is a finite one.
It is clear from the definition that $\mu_{\Op\ps{G}}(M) =
\mu_{\Op\ps{G}}(M(\pi))$. Also, it is not difficult to see that this
definition coincides with the classical notion of the
$\mu$-invariant for $\Ga$-modules when $G=\Ga\cong\Zp$. We now
record certain properties of this invariant which will be required
in the subsequent of the paper.

\bl \label{mu lemma} Let $G$ be a compact $p$-adic Lie group and let
$M$ be a finitely generated $\Op\ps{G}$-module. Then we have the
following statements.

\begin{enumerate}
\item[$(a)$] If $G_0$ is an open subgroup of $G$, then
 \[ [G:G_0]\mu_{\Op\ps{G}}(M) = \mu_{\Op\ps{G_0}}(M).\]

\item[$(b)$] Let $\Op'$ denote the ring of integers of a finite
extension $K'$ of $K$. Then we have
 \[ \mu_{\Op'\ps{G}}(M\ot_{\Op}\Op') = \mu_{\Op\ps{G}}(M).\]

\item[$(c)$] Viewing
$M$ as a $\Zp\ps{G}$-module, we have $[K:\Qp]\rank_{\Op\ps{G}}(M) =
\rank_{\Zp\ps{G}}(M)$.

\item[$(d)$] Viewing
$M$ as a $\Zp\ps{G}$-module, we have
$[k:\mathbb{F}_p]\mu_{\Op\ps{G}}(M) = \mu_{\Zp\ps{G}}(M)$.

\item[$(e)$] Suppose further that $G\cong \Zp^r$ for $r\geq 1$.
Then one has
\[\mu_{\Op\ps{G}}(M/\pi^n) =
n\rank_{\Op\ps{G}}(M) + \mu_{\Op\ps{G}}(M) \quad \mbox{for}~ n\gg 0.
\]
  \end{enumerate}
\el

\bpf (a), (b) and (c) are immediate from the definition. To prove
(d), it suffices, by (a), to consider the case when $G$ is a uniform
pro-$p$ group. By \cite[Corollary 1.7]{Ho}, one has
 \[\mu_{\Op\ps{G}}(M) = \displaystyle\sum_{i\geq 0}
(-1)^i\mathrm{ord}_q\big(H_i(G,M(\pi))\big),\] where $q$ is the
order of $k$. (d) is now immediate from this formula and the facts
that $M(\pi) = M(p)$ and that $[k:\mathbb{F}_p]\ord_q N = \ord_p N$.

To prove (e), we first consider the case when $M$ is a torsion
$\Op\ps{G}$-module. By the structure theory of torsion
$\Op\ps{G}$-module, we have that $M$ is pseudo-isomorphic to
\[ \bigoplus_{i=1}^s\Op\ps{G}/\pi^{\al_i}
\oplus \bigoplus_{j=1}^t\Op\ps{G}/f_j \] for some nonzero $f_j$
coprime to $\pi$. For $n \geq \max\{\al_1,...,\al_s\}$, the
$\Op\ps{G}$-module $M/\pi^n$ is pseudo-isomorphic to
\[  \bigoplus_{i=1}^s\Op\ps{G}/\pi^{\al_i}. \]
The equality of (e) in this case is immediate. Now we consider the
case when $M$ is a finitely generated $\Op\ps{G}$-module with
$\Op\ps{G}$-rank $r >0$. Denote $M_t$ to be the maximal torsion
$\Op\ps{G}$-submodule of $M$ and write $M_{tf} = M/M_t$. Since
$M_{tf}$ is torsionfree, we can find an injection $\Op\ps{G}^r
\hookrightarrow M_{tf}$ with a $\Op\ps{G}$-torsion cokernel $N$. By
the above argument, for sufficiently large $n$, one has
\[ \mu_{\Op\ps{G}}(M_t) = \mu_{\Op\ps{G}}(M_t/\pi^n).\] As the
$\mu_{\Op\ps{G}}$-invariant is additive on exact sequences of
torsion $\Op\ps{G}$-modules, the exact sequence
 \[ 0 \lra N[\pi^n]\lra N\stackrel{\pi^n}{\lra} N \lra N/\pi^n\lra 0\]
 shows that \[
\mu_{\Op\ps{G}}(N[\pi^n]) = \mu_{\Op\ps{G}}(N/\pi^n).\] Combining
the above equalities with the following two exact sequences
\[ 0\lra M_t/\pi^n \lra M/\pi^n
\lra M_{tf}/\pi^n\lra 0  \]
\[ 0\lra N[\pi^n] \lra \big(\Op\ps{G}/\pi^n\big)^r
\lra M_{tf}/\pi^n \lra N/\pi^n \lra 0,  \] we obtain the required
equality, noting that $\mu_{\Op\ps{G}}(M) = \mu_{\Op\ps{G}}(M_t)$.
\epf

Before stating the next lemma, we recall some terminology and
notation. Let $\rho: G\lra GL_d(\Op')$ be a continuous group
homomorphism with $\Op\subseteq \Op'$. Denote $W_{\rho}$ to be a
free $\Op'$-module of rank $d$ realizing $\rho$. If $M$ is a
$\Op\ps{G}$-module, we define $\mathrm{tw}_{\rho}(M)$ to be the
$\Op'$-module $W_{\rho}\ot_{\Zp}M$ with $G$ acting diagonally. We
shall say that $\rho$ is an \textit{Artin representation} if $\rho$
has finite image.

\bl \label{mu lemma Artin twist} Let $G$ be a compact $p$-adic Lie
group and let $M$ be a finitely generated $\Op\ps{G}$-module. Let
$\rho: G\lra GL_d(\Op')$ be an Artin representation with
$\Op\subseteq \Op'$. Then
 \[ \mu_{\Op'\ps{G}}(\mathrm{tw}_{\rho}(M)) = \mu_{\Op\ps{G}}(M)d.\]
\el

\bpf
 Let $G_0$ be a normal uniform pro-$p$ subgroup of $G$ such that
 $G_0 \subseteq \ker\rho$. By Lemma \ref{mu lemma}(a), it suffices
 to show that
 $\mu_{\Op'\ps{G_0}}(\mathrm{tw}_{\rho}(M)) =
 \mu_{\Op\ps{G_0}}(M)d$. By our choice of $G_0$, we have
 $\mathrm{tw}_{\rho}(M) = (M\ot_{\Op}\Op')^d$ as
 $\Op'\ps{G_0}$-modules. Thus, we have
 \[\mu_{\Op'\ps{G_0}}(\mathrm{tw}_{\rho}(M)) =
 \mu_{\Op'\ps{G_0}}(M\ot_{\Op}\Op') d =  \mu_{\Op\ps{G_0}}(M)d,  \]
 where the last equality follows from Lemma \ref{mu lemma}(b).
\epf

We now quote a result of Cucuo and Monsky \cite[Theorem 4.13]{CM}
(see also \cite[Theorem 3.12]{Mon} for a finer statement) which
extends the classical asymptotic formula to the case when $G\cong
\Zp^r$, $r\geq 2$.

\bt [Cucuo-Monsky] \label{Cucuo-Monsky}
 Suppose that $G\cong \Zp^r$, where $r\geq 2$.  Denote $G_m$ to be $G^{p^m}$.
Let $M$ be a finitely generated torsion $\Zp\ps{G}$-module which
satisfies the property that $\rank_{\Zp}M_{G_m} = O(p^{(r-2)m})$.
Then we have
 \[ \ord_p \big(M_{G_m}\big) = \mu_G(M)p^{rm} + l_0(M)mp^{(r-1)m}+ O(p^{(r-1)m})
 \quad \mbox{for}~m\gg 0, \]
 where $l_0(M)$ is defined as in \cite[Definition 1.2]{CM}.
\et

We will require a few more lemmas which estimate the order of
certain cohomology groups. For an abelian group $N$, we define its
\textit{$p$-rank} to be the $\mathbb{F}_p$-dimension of $N[p]$ which
we denote by $r_p(N)$. If $G$ is a pro-$p$ group, we write $h_1(G) =
r_p\big(H^1(G,\Z/p\Z)\big)$ and $h_2(G) =
r_p\big(H^2(G,\Z/p\Z)\big)$. We now state and prove the following
lemma which gives an estimate of the $p$-rank of the first and
second cohomology groups.

\bl \label{cohomology rank inequalities} Let $G$ be a pro-$p$ group,
and let $M$ be a discrete $G$-module which is cofinitely generated
over $\Zp$.
 If $h_1(G)$ is finite, then $r_p\big(H^1(G,M)\big)$ is finite, and
 we have the following estimate
 \[
 r_p\big(H^1(G,M)\big) \leq 2h_1(G) + \mathrm{corank}_{\Zp}(M) +
\ord_p( M/M_{\mathrm{div}}).
\]

If $h_2(G)$ is finite, then $r_p\big(H^2(G,M)\big)$ is finite, and
 we have the following estimate
 \[
 r_p\big(H^2(G,M)\big) \leq 2h_2(G) + \mathrm{corank}_{\Zp}(M) +
\ord_p( M/M_{\mathrm{div}}).
\] \el

\bpf We shall only give a proof of the upper bound for
$r_p\big(H^1(G,M)\big)$, the proof of the upper bound of
$r_p\big(H^2(G,M)\big)$ being similar. By \cite[Corollary
1.6.13]{NSW}, the only simple $G$-module is $\Z/p\Z$ with a trivial
$G$-action. If $M$ is finite, it follows from a standard d\'evissage
argument that
 \[ |H^1(G,M)| \leq |H^1(G,\Z/p\Z)|\ord_p(M).\]
 This in turn implies that
\[ r_p\big(H^1(G,M)\big) \leq \ord_p\big(H^1(G,M)\big) \leq h_1(G)+ \ord_p(M).
\]
  For a general $M$, we denote $M_{\mathrm{div}}$ to be the maximal $p$-divisible subgroup of
  $M$. Note that $M_{\mathrm{div}}$ is a $G$-submodule of $M$. Then
 we have a short exact sequence
\[ 0\lra M_{\mathrm{div}}\lra M \lra M/M_{\mathrm{div}}\lra 0\]
which induces the following exact sequence
\[ H^1(G, M_{\mathrm{div}})\lra
H^1(G, M) \lra H^1(G, M/M_{\mathrm{div}}).\]
  Therefore, we are
reduced to showing that $r_p\big(H^1(G, M_{\mathrm{div}})\big)$ and
$r_p\big(H^1(G, M/M_{\mathrm{div}})\big)$ are finite, and that the
following inequalities
\[ \ba{c}
r_p\big(H^1(G, M_{\mathrm{div}})\big)\leq h_1(G)\ + \mathrm{corank}_{\Zp}(M),\\
  r_p\big(H^1(G, M/M_{\mathrm{div}})\big)\leq
h_1(G) + \ord_p(M/M_{\mathrm{div}}) \ea \] hold. Since $M$ is
cofinitely generated over $\Zp$, we have that $ M/M_{\mathrm{div}}$
is finite. The finiteness of $r_p\big(H^1(G,
M/M_{\mathrm{div}})\big)$ and the validity of the second inequality
then follow from the above discussion. To see that the first
inequality holds, we first note that the short exact sequence
\[ 0\lra M_{\mathrm{div}}[p]\lra M_{\mathrm{div}} \stackrel{p}{\lra} M_{\mathrm{div}}\lra 0\] of
discrete $G$-modules induces a surjection $H^1(G,
M_{\mathrm{div}}[p]) \tha H^1(G, M_{\mathrm{div}})[p]$ and,
consequently, the inequality
$$ r_p\big( H^1(G, M_{\mathrm{div}})\big) \leq r_p\big(H^1(G,
M_{\mathrm{div}}[p])\big).
 $$
Again by the above discussion, the latter is less than or equal to
$h_1(G) +\ord_p(M_{\mathrm{div}}[p])$, and the required inequality
follows from the observation that $\ord_p(M_{\mathrm{div}}[p]) =
\mathrm{corank}_{\Zp}(M)$. \epf

As an immediate consequence of the preceding lemma, we have the
following.

\bl \label{cohomology order inequalities} Let $G$ be a pro-$p$
group, and let $M$ be a discrete $G$-module which is cofinitely
generated over $\Zp$.
 If $h_1(G)$ is finite, then
 we have the following estimate
 \[
 \ord_p\big(H^1(G,M)[p^n]\big) \leq n\big(2h_1(G) + \mathrm{corank}_{\Zp}(M) +
\ord_p( M/M_{\mathrm{div}})\big).
\]

If $h_2(G)$ is finite, then
 we have the following inequality
 \[
 \ord_p\big(H^2(G,M)[p^n]\big) \leq n\big(2h_2(G) + \mathrm{corank}_{\Zp}(M) +
\ord_p( M/M_{\mathrm{div}})\big).
\] \el

\section{A ratio formula} \label{a ratio formula}

As before, let $p$ be a prime, and let $F$ be a number field. If
$p=2$, we assume further that $F$ has no real primes. Suppose that
we are given the following data:

\begin{enumerate}
 \item[(a)] $M$ is a finite $\Gal(\bar{F}/F)$-module of $p$-power order
which is unramified outside a finite set of primes of $F$.

 \item[(b)] For each prime $v$ of $F$ above $p$, $M_v$ is a
$\Gal(\bar{F}_v/F_v)$-submodule of $M$.

\item[(c)] For each real prime $v$ of $F$, we write $M_v^+=
M^{\Gal(\bar{F}_v/F_v)}$.

\item[(d)] The following equality
 \begin{equation} \label{finite equality}
 |M|^{r_2(F)}\displaystyle\prod_{v~\mathrm{real}}|M/M_v^+| = \prod_{v|p}|M/M_v|
 \end{equation}
holds. Here $r_2(F)$ denotes the number of complex primes of $F$.
\end{enumerate}

We will denote the data by $\big(M, \{M_v\}_{v|p},
\{M_v^+\}_{v|\R}\big)$. Let $S$ denote a finite set of primes of $F$
which contains all the primes above $p$, the ramified primes of $M$
and all infinite primes. Denote $F_S$ to be the maximal algebraic
extension of $F$ unramified outside $S$. For each algebraic
extension $\mathcal{L}$ of $F$ contained in $F_S$, we write
$G_S(\mathcal{L}) = \Gal(F_S/\mathcal{F})$.

\bl \label{euler char equality} For a given data $\big(M,
\{M_v\}_{v|p}, \{M_v^+\}_{v|\R}\big)$, we have the following
equality
 \[ \frac{|H^0(G_S(F),M)| |H^2(G_S(F),M)|}{|H^1(G_S(F),M)|}
 = \prod_{v|p} \frac{|H^0(F_v,M/M_v)| |H^2(F_v,M/M_v)|}{|H^1(F_v,M/M_v)|}.  \]
\el

\bpf
 This follows from a straightforward application of the local and global Euler
characteristic formulas. \epf

Now set
\[ H^1_{str}(F_v, M)=
\begin{cases} \ker\big(H^1(F_v, M)\lra H^1(F_v, M/M_v)\big) & \text{\mbox{if} $v|p$},\\
 \ker\big(H^1(F_v, M)\lra H^1(F^{ur}_v, M)\big) & \text{\mbox{if} $v\nmid p$,}
\end{cases} \]
 where $F_v^{ur}$ is the maximal unramified extension of $F_v$.
The Greenberg strict Selmer group (see \cite[P. 116]{G89}) attached
to  $\big(M, \{M_v\}_{v|p}, \{M_v^+\}_{v|\R}\big)$ is defined by
\[ S(M/F):= \Sel^{str}(M/F) := \ker\Big( H^1(G_S(F),M)\lra \bigoplus_{v \in S}H^1(F_v,
M)/H^1_{str}(F_v, M)\Big).\]

For the remainder of the paper, we shall refer the Greenberg strict
Selmer group as the Selmer group. For a data  $\big(M,
\{M_v\}_{v|p}, \{M_v^+\}_{v|\R}\big)$, we say that $\big(M^*,
\{M^*_v\}_{v|p}, \{(M^*)_v^+\}_{v|\R}\big)$ is the dual data given
by $M^* = \Hom_{\cts}(M,\mu_{p^{\infty}})$, $M_v^*=
\Hom_{\cts}(M/M_v,\mu_{p^{\infty}})$ and  $(M^*)_v^+=
\Hom_{\cts}(M/M^+_v,\mu_{p^{\infty}})$. It is a straightforward
exercise to verify that the dual data satisfies (\ref{finite
equality}).  The Selmer group $S(M^*/F)$ for  $\big(M^*,
\{M^*_v\}_{v|p}, \{(M^*)_v^+\}_{v|\R}\big)$ is defined similarly.
Then one has the following proposition.

\bp \label{ratio of Selmer}
 For a data  $\big(M,
\{M_v\}_{v|p}, \{M_v^+\}_{v|\R}\big)$, we have the following
equality
 \[ \frac{|S(M/F)|
 \prod_{v|p}|H^0(F_v,M/M_v)|}{|H^0(G_S(F),M)|}
 = \frac{|S(M^*/F)| \prod_{v|p}|H^0(F_v,M^*/M^*_v)|}{ |H^0(G_S(F),M^*)|}.  \]
\ep

\bpf
 This is proven in the same way as \cite[Formula
 (53)]{G89} by combining a Poitou-Tate sequence argument with Lemma \ref{euler char equality}.
\epf

\section{Comparing Selmer groups over multiple $\Zp$-extensions}
\label{compare}

As before, let $p$ be a prime. We let $F$ be a number field. If
$p=2$, we assume further that $F$ has no real primes. Denote $\Op$
to be the ring of integers of some finite extension $K$ of $\Qp$.
Fix a local parameter $\pi$ for $\Op$. Suppose that we are given the
following data:

\begin{enumerate}
 \item[(a)] $A$ is a
cofinitely generated cofree $\Op$-module of $\Op$-corank $d$ with a
continuous, $\Op$-linear $\Gal(\bar{F}/F)$-action which is
unramified outside a finite set of primes of $F$.

 \item[(b)] For each prime $v$ of $F$ above $p$, $A_v$ is a
$\Gal(\bar{F}_v/F_v)$-submodule of $A$ which is cofree of
$\Op$-corank $d_v$

\item[(c)] For each real prime $v$ of $F$, we write $A_v^+=
A^{\Gal(\bar{F}_v/F_v)}$.

\item[(d)] The following equality
\begin{equation} \label{data equality}
  \sum_{v|p} (d-d_v)[F_v:\Qp] = dr_2(F) +
 \sum_{v~\mathrm{real}}(d-d^+_v)
  \end{equation}
holds. Here $r_2(F)$ denotes the number of complex primes of $F$.
\end{enumerate}

We denote the above data as $\big(A, \{A_v\}_{v|p}, \{A^+_v\}_{v|\R}
\big)$. From these data, we define its dual data as follows. Set
$A^* = \Hom_{\cts}(T_{\pi}(A),\mu_{p^{\infty}})$, $A^*_v=
\Hom_{\cts}(T_{\pi}(A/A_v),\mu_{p^{\infty}})$ and $(A^*)^+_v=
\Hom_{\cts}(T_{\pi}(A/A^+_v),\mu_{p^{\infty}})$. Here $T_{\pi}(N)$
denotes the $\pi$-adic Tate module of a $\Op$-module $N$. It is an
easy exercise to verify that $\big(A, \{A_v\}_{v|p},
\{(A^*)^+_v\}_{v|\R} \big)$ satisfies equality (\ref{data
equality}). For each $n$ and an $\Op$-module $N$, we denote
$N[\pi^n]$ to be the kernel of $\pi^n:N\lra N$. One can check easily
that the induced data $\big(A[\pi^n], \{A_v[\pi^n]\}_{v|p},
\{A^+_v[\pi^n]\}_{v|\R}\big)$ satisfies equality (\ref{finite
equality}).

We now describe briefly the general arithmetic situation, where we
can obtain the above data from. (See Section \ref{examples} for more
explicit examples.) Let $V$ be a $d$-dimensional $K$-vector space
with a continuous $G_S(F)$-action. Suppose that for each prime $v$
of $F$ above $p$, there is a $d_v$-dimensional $K$-subspace $V_v$ of
$V$ which is invariant under the action of $\Gal(\bar{F}_v/F_v)$,
and for each real prime $v$ of $F$, $V^{\Gal(\bar{F}_v/F_v)}$ has
dimension $d_v^+$. Choose a $G_S(F)$-stable $\Op$-lattice $T$ of $V$
(Such a lattice exists by compactness). We can obtain a data from
$V$ by setting $A= V/T$ and $A_v = V_v/ T\cap V_v$. Note that $A$
and $A_v$ depends on the choice of the lattice $T$.

We now consider the base change property of our data. Let $L$ be a
finite extension of $F$. We obtain another data $\big(A,
\{A_w\}_{w|p}, \{A^+_w\}_{w|\R} \big)$ over $L$ as follows: we
consider $A$ as a $\Gal(\bar{F}/L)$-module, and for each prime $w$
of $L$ above $p$, we set $A_w =A_v$, where $v$ is a prime of $F$
below $w$, and view it as a $\Gal(\bar{F}_v/L_w)$-module. Then $d_w
= d_v$. For each real prime $w$ of $L$, one sets
$A^{\Gal(\bar{L}_w/L_w)}= A^{\Gal(\bar{F}_v/F_v)}$ and writes $d^+_w
= d^+_v$, where $v$ is a real prime of $F$ below $w$. In general,
the $d_w$'s and $d_w^+$ need not satisfy equality (\ref{data
equality}). We now record the following lemma which give sufficient
condition for equality (\ref{data equality}) to hold for the data
$\big(A, \{A_w\}_{w|p}, \{A^+_w\}_{w|\R} \big)$ over $L$.

\bl \label{data base change} Suppose that $\big(A, \{A_v\}_{v|p},
\{A^+_v\}_{v|\R} \big)$ is a data defined over $F$.
 Suppose further that at least one of the following statements holds.
 \begin{enumerate}
\item[$(i)$] All the archimedean primes of $F$ are unramified in
$L$.

\item[$(ii)$] $[L:F]$ is odd

\item[$(iii)$] $F$ is totally imaginary.

\item[$(iv)$] $F$ is totally real, $L$ is totally imaginary and
\[ \sum_{v~\mathrm{real}} d^+_v = d[F:Q]/2.\]
 Then we have the equality
 \[ \sum_{w|p} (d-d_w)[L_w:\Qp] = dr_2(L) +
 \sum_{w~\mathrm{real}}(d-d^+_w).\]
\end{enumerate}
 \el

\bpf Note that if either of the assertions in (ii) or (iii) holds,
then the assertion in (i) holds. Therefore, to prove the lemma in
these cases, it suffices to prove it under the assumption of (i). We
have the following calculation
 \[ \ba{rl}
\displaystyle \sum_{w|p} (d-d_w)[L_w:\Qp]\!\!
   &= \displaystyle\sum_{v|p}\sum_{w|v} (d-d_v)[L_w:F_v][F_v:\Qp] \\
   &= \displaystyle\sum_{v|p} (d-d_v)[F_v:\Qp] \sum_{w|v}[L_w:F_v] \\
   &= \displaystyle [L:F]\sum_{v|p} (d-d_v)[F_v:\Qp] \\
  &= \displaystyle d[L:F]r_2(F) + [L:F]\sum_{v~\mathrm{real}} (d-d^+_v). \\  \ea \]
 Since (i) holds, every prime of $L$ above a real prime (resp.,
complex prime) of $F$ is a real prime (resp., complex prime).
Therefore, one has $[L:F]r_2(F) = r_2(L)$ and
\[ [L:F]\sum_{v~\mathrm{real}} (d-d^+_v) = \sum_{w~\mathrm{real}} (d-d^+_w). \]
The required conclusion then follows.

Now suppose that (iv) holds. Then $r_2(F) =0$ and the above sum is
 \[ \ba{rl}
\displaystyle [L:F]\sum_{v~\mathrm{real}} (d-d^+_v)\!\! &=
  \displaystyle[L:F]\sum_{v~\mathrm{real}}d - [L:F]\sum_{v~\mathrm{real}}d^+_v \\
   &= [L:\Q]d - [L:F] d[F:\Q]/2 \vspace{0.1in} \\
 &= d[L:\Q]/2 = dr_2(L).
\ea
\]
\epf

Let $S$ be a finite set of primes of $F$ which contains all the
primes above $p$, the ramified primes of $A$ and all infinite
primes. Denote $F_S$ to be the maximal algebraic extension of $F$
unramified outside $S$ and write $G_S(\mathcal{L}) =
\Gal(F_S/\mathcal{L})$ for every algebraic extension $\mathcal{L}$
of $F$ which is contained in $F_S$. Let $L$ be a finite extension of
$F$ contained in $F_S$ such that the data $\big(A, \{A_w\}_{w|p},
\{A^+_w\}_{w|\R} \big)$ satisfies (\ref{data equality}). For a prime
$w$ of $L$ lying over $S$, set
\[ H^1_{str}(L_w, A)=
\begin{cases} \ker\big(H^1(L_w, A)\lra H^1(L_w, A/A_w)\big) & \text{\mbox{if} $w$
 divides $p$},\\
 \ker\big(H^1(L_w, A)\lra H^1(L^{ur}_w, A)\big) & \text{\mbox{if} $w$ does not divide $p$,}
\end{cases} \]
 where $L_w^{ur}$ is the maximal unramified extension of $L_w$.
The (Greenberg strict) Selmer group attached to the data is then
defined by
\[ S(A/L) := \Sel^{str}(A/L) := \ker\Big( H^1(G_S(L),A)\lra
\bigoplus_{w \in S_L}H^1_s(L_w, A)\Big),\] where we write
$H^1_s(L_w, A) = H^1(L_w, A)/H^1_{str}(L_w, A)$ and $S_L$ denotes
the set of primes of $L$ above $S$. It is straightforward to verify
that $S(A/L) = \ilim_n S(A[\pi^n]/L)$, where $S(A[\pi^n]/L)$ the
Selmer group of $\big(A[\pi^n], \{A_w[\pi^n]\}_{w|p},
\{A^+_w[\pi^n]\}_{w|\R} \big)$ defined in Section \ref{a ratio
formula}, and the direct limit is taken over the maps $S(A[\pi^n]/L)
\lra S(A[\pi^{n+1}]/L)$ induced by the natural injections
$A[\pi^n]\hookrightarrow A[\pi^{n+1}]$ and
$A_w[\pi^n]\hookrightarrow A_w[\pi^{n+1}]$. We will write $X(A/L)$
for its Pontryagin dual. The $S(A^*/L)$ are defined analogously and
one also has the identification $S(A^*/L) = \ilim_n
S(A^*[\pi^n]/L)$.

\smallskip
For an infinite algebraic extension $\mathcal{L}$ of $F$ contained
in $F_S$, we define $S(A/\mathcal{L}) = \ilim_L S(A/L)$, where the
limit runs over all finite extensions $L$ of $F$ contained in
$\mathcal{L}$. We will write $X(A/\mathcal{L})$ for the Pontryagin
dual of $S(A/\mathcal{L})$.

\smallskip
Let $F_{\infty}$ be a $\Zp^r$-extension of $F$ which contains the
cyclotomic $\Zp$-extension $F_{\cyc}$ of $F$. We denote $G$ to be
the Galois group. Let $F_m$ be the unique subextension of
$F_{\infty}$ over $F$ with $\Gal(F_{m}/F) \cong (\Z/p^{m})^r$ and
write $G_m =\Gal(F_{\infty}/F_m)$. We can now prove the following
theorem which is the main result of this section.

\bt \label{Zp^r} Retain the above assumptions.
 $X(A/F_{\infty})$ and $X(A^*/F_{\infty})$ have the same
 $\Op\ps{G}$-rank and the same $\mu_{\Op\ps{G}}$-invariant.
\et

\bpf
 The proof follows the argument in \cite{G89}. As the case of $r=1$ is essentially
dealt there, we will concentrate on the case when $r\geq 2$. By
Lemma \ref{mu lemma}(c)(d), it suffices to show the theorem under
the assumption that $\Op=\Zp$. Fix an arbitrary positive integer
$n$. By Lemma \ref{data base change}(ii) and (iii), for each $m$, we
may apply Proposition \ref{ratio of Selmer} to obtain the equality
\[ \frac{|S(A[p^n]/F_m)|
 \prod_{v_m|p}|H^0(F_{m, v_m},A[p^n]/A_{v_m}[p^n])|}{|H^0(G_S(F_m),A[p^n])|}
 = \frac{|S(A^*[p^n]/F_m)| \prod_{v_m|p}|H^0(F_{m,v_m}, A^*[p^n]/A_{v_m}^*[p^n])|}
 { |H^0(G_S(F_m),A^*[p^n])|},  \]
 where $v_m$ runs over all the primes of $F_m$ above $p$.
Clearly, $|H^0(G_S(F_m),A[p^n])|$ and $|H^0(G_S(F_m),A^*[p^n])|$ are
bounded independent of $m$ (for a fixed $n$). Since there are only
finite number of primes of $F_{\cyc}$ above $p$, the decomposition
group of $v$ in $G$ has at most dimension $r-1$. Therefore, it
follows that $\prod_{v_m|p}|H^0(F_{m,v_m}, A[p^n]/A_v[p^n])|$ and
$\prod_{v_m|p}|H^0(F_{m,v_m}, A^*[p^n]/A_v^*[p^n])|$ are both
$p^{O(p^{(r-1)m})}$. Thus, we have
 \begin{equation} \label{A and A*} \ord_p \big(S(A[p^n]/F_m)\big) = \ord_p\big(S(A^*[p^n]/F_m)\big)
+O(p^{(r-1)m}). \end{equation}
 Now we need to estimate the order of the kernels and
cokernels of the maps
 \[ S(A[p^n]/F_m) \stackrel{r_m}{\lra} S(A/F_m)[p^n]
 \stackrel{s_m}{\lra}  \big(S(A/F_{\infty})[p^n]\big)^{G_m}.\]
 One sees easily that $\ker r_m \subseteq A(F_m)/p^n$ and $\ker s_m
 \subseteq H^1(G_m, A(F_{\infty}))[p^n]$. It is clear that one has $\ord_p(\ker r_m) =
 O(1)$. On the other hand, it follows from Lemma \ref{cohomology order inequalities}
that $\ord_p\big(H^1(G_m, A(F_{\infty}))[p^n]\big) = O(1)$ (noting
that $h_1(G_m)$ is a constant function in $m$). Thus, one has
$\ord_p(\ker s_m) =  O(1)$.

To estimate $\coker r_m$ and $\coker s_m$, one first observes that
$\ord_p(\coker r_m) \leq \ord_p(\ker r_m')$ and that $\ord_p(\coker
s_m) \leq \ord_p(\ker s_m') + \ord_p(H^2(G_m, A(F_{\infty}))[p^n])$,
where $r_m'$ and $s_m'$ are given by
\[r_m' = \big(r'_{m,v_m}\big): \bigoplus_{v_m\in S(F_m)}H^1_s(F_{m,v_m}, A[p^n])
\lra \bigoplus_{v_m\in S(F_m)}H^1_s(F_{m,v_m}, A)[p^n];\]
\[s_m'  = \big(s'_{m,v_m}\big) : \bigoplus_{v_m\in S(F_m)}H^1_s(F_{m,v_m}, A)[p^n]
\lra \Big(\ilim_m\bigoplus_{v_m\in S(F_m)}H^1_s(F_{m,v_m},
A)[p^n]\Big)^{G_m}.\]
 By Lemma \ref{cohomology order inequalities}, one has that
$\ord_p\big(H^2(G_m, A(F_{\infty}))[p^n]\big) = O(1)$ (noting that
$h_2(G_m)$ is a constant function in $m$). To estimate $\coker
r'_m$, we first observe that

\[\ker r'_{m,v_m} \subseteq
\begin{cases} \ker\big(H^1(F_{m, v_m}, A/A_{v_m}[p^n])\lra
 H^1(F_{m, v_m}, A/A_{v_m})[p^n]
 & \text{\mbox{if} $v_m|p$},\\
\ker \big(H^1(F^{ur}_{m, v_m}, A[p^n])\lra
 H^1(F^{ur}_{m, v_m}, A)[p^n]  & \text{\mbox{if} $v_m\nmid p$,}
\end{cases} \]

\[  =
\begin{cases} A/A_{v_m}(F_{m, v_m})[p^n] \hspace{0.6in}
 & \text{\mbox{if} $v_m|p$}, \hspace{1in}\\
 A(F^{ur}_{m, v_m})[p^n] \hspace{0.6in}
 & \text{\mbox{if} $v_m\nmid p$, \hspace{1in}}
\end{cases} \]

It is now clear from the above that $\ord_p(\ker r'_{m,v_m})$ is
bounded independent of $m$ and $v_m$ (for a fixed $n$). Combining
these estimates with the fact that the decomposition group of $v$ in
$G$ has dimension at most $r-1$ for every $v\in S$, one then has the
estimate $\ord_p \big(\ker r'_m\big) = O(p^{(r-1)m})$.

 To estimate $\coker
s'_m$, we now observe that

\[\ker s'_{m,v_m} \subseteq
\begin{cases} H^1\big(G_{m, v_m}, A/A_{v_m}(F_{m,v_m})\big)[p^n]
 & \text{\mbox{if} $v_m|p$},\\
H^1\big(\Gal(F_{\infty, v_m}/F_{m,v_m}^{ur}), A(F_{m,v_m})\big)[p^n]
& \text{\mbox{if} $v_m\nmid p$.}
\end{cases} \]

By Lemma \ref{cohomology order inequalities}, one can verify that
$\ker s'_{m,v_m}$ is bounded independent of $m$ and $v_m$ (for a
fixed $n$). As before, combining these estimates with the fact that
the decomposition group of $v$ in $G$ has dimension at most $r-1$
for every $v\in S$, we obtain
  $\ord_p \big(\ker s'_m\big)
= O(p^{(r-1)m})$. In conclusion, we have
 \begin{equation} \label{A} \ord_{p} \big(S(A[p^n]/F_m)\big) =
 \ord_{p}\big(S(A/F_{\infty})[p^n]^{G_m}\big) +O(p^{(r-1)m}).
  \end{equation}
Similarly, one also has
 \begin{equation} \label{A*} \ord_{p} \big(S(A^*[p^n]/F_m)\big) =
 \ord_{p}\big(S(A^*/F_{\infty})[p^n]^{G_m}\big) +O(p^{(r-1)m}).
  \end{equation}
Combining the estimates in (\ref{A and A*}), (\ref{A}) and
(\ref{A*}), we obtain
 \begin{equation} \label{/p} \ord_p\Bigg(\Big(X(A/F_{\infty})/p^n\Big)_{G_m}\Bigg) =
\ord_p\Bigg(\Big(X(A^*/F_{\infty})/p^n\Big)_{G_m}\Bigg)
+O(p^{(r-1)m}).\end{equation}
 Since $X(A/F_{\infty})/p^n$ and $X(A^*/F_{\infty})/p^n$
are $p$-torsion modules, we have
\[ \rank_{\Zp}\Big(X(A/F_{\infty})/p^n\Big)_{G_m} = \rank_{\Zp}
\Big(X(A^*/F_{\infty})/p^n\Big)_{G_m} =0,\] and therefore, the
hypothesis of Theorem \ref{Cucuo-Monsky} is satisfied. Hence we may
combine the said theorem with the estimate in (\ref{/p}) to conclude
that
\[
 \mu_{\Zp\ps{G}}\Big(X(A/F_{\infty})/p^n\Big)p^{rm}=
\mu_{\Zp\ps{G}}\Big(X(A^*/F_{\infty})/p^n\Big)p^{rm} + O(p^{(r-1)m})
\] (note that since $X(A/F_{\infty})/p^n$ and $X(A^*/F_{\infty})/p^n$ are $p$-torsion modules, the
$l_0$ quantity attached to these modules vanishes) which in turn
implies the equality
\[ \mu_{\Zp\ps{G}}\Big(X(A/F_{\infty})/p^n\Big) =
\mu_{\Zp\ps{G}}\Big(X(A^*/F_{\infty})/p^n\Big).\] By Lemma \ref{mu
lemma}(e), this in turn implies that
\[ n\rank_{\Zp\ps{G}}\big(X(A/F_{\infty})\big) + \mu_{\Zp\ps{G}}\Big(X(A/F_{\infty})\Big) =
n\rank_{\Zp\ps{G}}\big(X(A^*/F_{\infty})\big)+
\mu_{\Zp\ps{G}}\Big(X(A^*/F_{\infty})\Big)\] for $n\gg 0$.
Therefore, we have the equalities $\rank_{\Zp\ps{G}}
\Big(X(A/F_{\infty})\Big)=
\rank_{\Zp\ps{G}}\Big(X(A^*/F_{\infty})\Big)$ and
$\mu_{\Zp\ps{G}}\Big(X(A/F_{\infty})\Big)=
\mu_{\Zp\ps{G}}\Big(X(A^*/F_{\infty})\Big)$. \epf

We record several immediate corollaries.

\bc \label{Zp^r corollary}
 $X(A/F_{\infty})$ is a $\Op\ps{G}$-torsion module if and only
 if $X(A^*/F_{\infty})$ is a $\Op\ps{G}$-torsion module.  \ec

\bc \label{Zp corollary}
 Suppose that $G=\Ga\cong \Zp$. Then
 $X(A/F_{\cyc})$ is a finitely generated $\Op$ module if and only
 if $X(A^*/F_{\cyc})$ is a finitely generated $\Op$ module.  \ec

\br
 For most data coming from (nearly) ordinary
representations, it is expected that $X(A/F_{\cyc})$ is a torsion
$\Op\ps{\Ga}$-module (see \cite[Conjecture 1]{G89} or
\cite[Conjecture 1.7]{We}). When the data $(A, \{A_v\}_{v|p})$ comes
from the Galois representation attached to a primitive Hecke
eigenform for $GL_2$ over $\Q$ and when the base field $F$ is
abelian over $\Q$, the torsionness condition (over $F_{\cyc}$) is a
deep theorem of Kato \cite{Ka}. \er

\br
 If $F_{\infty}$ is a general $\Zp^r$-extension of $F$
(that does not contain $F_{\cyc}$) which has the property such that
for each prime $v\in S$, the decomposition group of
$\Gal(F_{\infty}/F)$ at $v$ has dimension $\leq r-1$, then the
argument of Theorem \ref{Zp^r} carries over to yield the same
conclusion. \er

\section{Comparing Selmer groups over noncommutative $p$-adic Lie extensions}
\label{compare2}

In this section, we will compare the Selmer groups of $A$ and $A^*$
over a noncommutative $p$-adic Lie extensions. We shall say that
$F_{\infty}$ is an \textit{admissible $p$-adic Lie extension} of $F$
if (i) $\Gal(F_{\infty}/F)$ is compact $p$-adic Lie group, (ii)
$F_{\infty}$ contains the cyclotomic $\Zp$ extension $F^{\cyc}$ of
$F$ and (iii) $F_{\infty}$ is unramified outside a set of finite
primes. Write $G = \Gal(F_{\infty}/F)$, $H =
\Gal(F_{\infty}/F_{\cyc})$ and $\Ga =\Gal(F_{\cyc}/F)$. Let $S$
denote a finite set of primes of $F$ which contains all the primes
above $p$, the ramified primes of $A$, the infinite primes and the
primes that are ramified in $F_{\infty}/F$.

Let $\big(A, \{A_v\}_{v|p}, \{A^+_v\}_{v|\R}\big)$ denote the data
defined in Section \ref{compare} and $\big(A^*, \{A^*_v\}_{v|p},
\{(A^*)^+_v\}_{v|\R}\big)$ the dual data. The Selmer group of $A$
over $F_{\infty}$ is defined to be $S(A/F_{\infty}) = \ilim_L
S(A/L)$, where $L$ runs through all finite extensions of $F$
contained in $F_{\infty}$. As before, we denote the Pontryagin dual
of $S(A/F_{\infty})$ by $X(A/F_{\infty})$. We have similar
definitions for $S(A^*/F_{\infty})$ and $X(A^*/F_{\infty})$.

\textbf{We shall impose the following hypothesis on our data
$\big(A, \{A_v\}_{v|p}, \{A^+_v\}_{v|\R}\big)$ and $p$-adic Lie
extension $F_{\infty}$ for the rest of this section: for every
finite extension $L$ of $F$ contained in $F_{\infty}$, the induced
data $\big(A, \{A_w\}_{w|p}, \{A^+_w\}_{w|\R}\big)$ over $L$ also
satisfies equation (\ref{data equality}).}

As a start, we record the following special case.

\bp \label{fg H}
 $X(A/F_{\infty})$ is a finitely generated
$\Op\ps{H}$-module if and only if $X(A^*/F_{\infty})$ is a finitely
generated $\Op\ps{H}$-module.
 \ep

\bpf
 Let $L$ be a finite extension of $F$ contained in $F_{\infty}$ such
that $\Gal(F_{\infty}/L)$ is pro-$p$. Since
$H_L=\Gal(F_{\infty}/L_{\cyc})$ is a subgroup of $H$ of finite
index, a $\Op\ps{G}$-module is finitely generated over $\Op\ps{H}$
if and only if it is finitely generated $\Op\ps{H_L}$-module. Now by
a standard argument (for instance, see Lemma \cite[Lemma
2.4]{CS12}), the natural map $X(A/F_{\infty})_{H_L}\lra
X(A/L_{\cyc})$ has kernel and cokernel that are finitely generated
over $\Op$. Since $H_L$ is pro-$p$, one can apply a Nakayama lemma
argument to conclude that $X(A/F_{\infty})$ is finitely generated
over $\Op\ps{H_{L}}$ if and only if $X(A/L_{\cyc})$ is finitely
generated over $\Op$. The conclusion of the proposition now follows
from the above discussion and Corollary \ref{Zp corollary}. \epf

We say that the admissible extension $F_{\infty}$ of $F$ is
\textit{almost abelian} if there exists a finite extension $L$ of
$F$ contained in $F_{\infty}$ such that $\Gal(F_{\infty}/L)\cong
\Zp^r$ for some $r\geq 1$. We then have the following theorem which
will follow from Theorem \ref{Zp^r} and our extended definition of
rank and $\mu$-variant.

\bt \label{almost abelian} Suppose that $F_{\infty}$ is an almost
abelian $S$-admissible extension of $F$. Then $X(A/F_{\infty})$ and
$X(A^*/F_{\infty})$ have the same $\Op\ps{G}$-rank and the same
$\mu_{\Op\ps{G}}$-invariant. In particular,  $X(A/F_{\infty})$ is a
$\Op\ps{G}$-torsion module if and only
 if $X(A^*/F_{\infty})$ is a $\Op\ps{G}$-torsion module.
 \et

To obtain a similar conclusion for the $\mu$-invariant over a
general noncommutative $p$-adic Lie extension, we need to assume a
stronger condition which was first introduced in \cite{CFKSV} and
was crucial in the formulation of the main conjectures of
non-commutative Iwasawa theory (see \cite{CFKSV, FK}). For a
finitely generated torsion $\Op\ps{G}$-module $M$, we say that $M$
\textit{belongs to} $\M_H(G)$ if $M/M(\pi)$ is finitely generated
over $\Op\ps{H}$. Here $M(\pi)$ is the submodule of $M$ consisting
of elements of $M$ annihilated by a power of $\pi$.

By the definition of the Selmer group, we have an exact sequence
\[ 0\lra S(A/F_{\infty}) \lra H^1(G_S(F_{\infty}),A) \stackrel{\la_{A/F_{\infty}}}{\lra}
 \bigoplus_{v \in S} J_v(A/F_{\infty}), \]
 where $J_v(A/F_{\infty}) =
 \ilim_L\bigoplus_{w|v} H^1_s(L_w, A)$.

\bt \label{MHG}  Let $F_{\infty}$ be an $S$-admissible $p$-adic Lie
extension. Assume that both $X(A/F_{\infty})$ and
$X(A^*/F_{\infty})$ belong to $\M_H(G)$.  Furthermore, suppose that
$A(L_{\cyc})$ and $A^*(L_{\cyc})$ are finite for
 every finite extension $L$ of $F$ contained in
$F_{\infty}$.

Then $X(A/F_{\infty})$ and $X(A^*/F_{\infty})$ have the same
$\mu_{\Op\ps{G}}$-invariant. \et



To prove Theorem \ref{MHG}, we require two lemmas.

\bl \label{torsion implies surjective} Let $F_{\infty}/F$ be an
admissible $p$-adic Lie extension. Assume that $X(A/F_{\infty})$
belongs to $\M_H(G)$. Furthermore, suppose that $A^*(L_{\cyc})$ is
finite for
 every finite extension $L$ of $F$ contained in
$F_{\infty}$. Then $H^2(G_S(F_{\infty}), A) =0$ and
$\la_{A/F_{\infty}}$ is surjective. \el

\bpf Since $X(A/F_{\infty})$ belongs to $\M_H(G)$, it follows from
the argument of \cite[Proposition 2.5]{CS12} that $X(A/L_{\cyc})$ is
a torsion $\Op\ps{\Ga_L}$-module for every finite extension $L$ of
$F$ contained in $F_{\infty}$, where $\Ga_L = \Gal(L_{\cyc}/L)$. The
Cassel-Poitou-Tate sequence gives an exact sequence
 \[ \ba{c}
  0\lra S(A/L_{\cyc})\lra H^1(G_S(L_{\cyc}), A)
 \stackrel{\la_{A/L_{\cyc}}}{\lra}
 \displaystyle\bigoplus_{w\in S(L_{\cyc})}H^1_s(L_w, A)
   \hspace{1in} \\
  \hspace{1in} \lra
 \big(\widehat{S}(A^*/L_{\cyc})\big)^{\vee}
 \lra H^2(G_S(L_{\cyc}), A)\lra \displaystyle\bigoplus_{w\in S(L_{\cyc})}H^2(L_w,
 A).
 \ea \]
  Here $\widehat{S}(A^*/L_{\cyc})$ is defined as the kernel of the
map
 \[ \plim_{E} H^1(G_S(L), T_{\pi} A^*) \lra \plim_{E}
 \bigoplus_{w|S}T_{\pi} H^1(L_w , A^*), \] where the inverse limit is
taken over all finite extensions $E$ of $L$ contained in $L_{\cyc}$.
Noting that $\displaystyle\bigoplus_{w\in S(L_{\cyc})}H^2(L_w, A)$
is cofinitely generated over $\Op$, it then follows from a
straightforward $\Op\ps{\Ga_L}$-rank calculation that
$\widehat{S}(A^*/L_{\cyc})$ has zero $\Op\ps{\Ga_L}$-rank. On the
other hand, by a similar argument to that in \cite[Proposition
7.1]{HV}, one has an injection
\[  \entrymodifiers={!! <0pt, .8ex>+} \SelectTips{eu}{}\xymatrix{
\widehat{S}(A^*/L_{\cyc})\big)~ \ar@{^{(}->}[r] &
\Hom_{\Op\ps{\Gal(L_{\cyc}/L)}}\big(S(A/L_{\cyc}),
\Op\ps{\Gal(L_{\cyc}/L)}\big)} \] which in turn implies that
$\widehat{S}(A^*/L_{\cyc})$ is torsionfree over $\Op\ps{\Ga_L}$.
Hence we must have $\widehat{S}(A^*/L_{\cyc}) =0$ and this in turn
implies that $\la_{A/L_{\cyc}}$ is surjective. Since
$\displaystyle\bigoplus_{w\in S(L_{\cyc})}H^2(L_w, A)$ is cofinitely
generated over $\Op$, it follows from the surjectivity of
$\la_{A/L_{\cyc}}$ and the above exact sequence that
$H^2(G_S(L_{\cyc}),A)$ is cotorsion over $\Op\ps{\Ga_L}$. On the
other hand, by \cite[Proposition 4]{G89}, $H^2(G_S(L_{\cyc}),A)$ is
a cofree $\Op\ps{\Ga_L}$-module. Therefore, this will force
$H^2(G_S(L_{\cyc}), A) = 0$. Since $H^2(G_S(F_{\infty}), A) =\ilim_L
H^2(G_S(L^{\cyc}), A)$ and $\la_{A/F_{\infty}} = \ilim_L
\la_{A/L_{\cyc}}$, it follows that $H^2(G_S(F_{\infty}), A) = 0$ and
$\la_{A/F_{\infty}}$ is surjective. \epf

\bl \label{MHG lemma} Suppose that $G$ is a compact pro-$p$ $p$-adic
Lie group with no $p$-torsion. Assume that $X(A/F_{\infty})$ belongs
to $\M_H(G)$. Suppose also that $H^2(G_S(F_{\cyc}), A) = 0$,
$H^2(G_S(F_{\infty}), A) = 0$, and $\la_{A/F_{\cyc}}$ and
$\la_{A/F_{\infty}}$ are surjective.

Then $\mu_{\Op\ps{G}}\big(X(A/F_{\infty})\big) =
\mu_{\Op\ps{\Ga}}\big(X(A/F_{\cyc})\big)$. \el

\bpf
 Under the assumptions of the lemma, one can apply a similar
argument to that of \cite[Proposition 2.13]{CSS} (see also
\cite[Proposition 4.7]{Lim2}) to conclude that
\[ \mu_{\Op\ps{G}}\big(X(A/F_{\infty})\big) = \mu_{\Op\ps{\Ga}}\big(X(A/F_{\cyc})\big) +
\sum_{n\geq 0}(-1)^{n+1}
\mu_{\Op\ps{\Ga}}\big(H_n(H,X_f(A/F_{\infty}))\big),
\]
 where $X_f(A/F_{\infty})=
X(A/F_{\infty})/X(A/F_{\infty})(\pi)$. By the hypothesis that
$X(A/F_{\infty})$ belongs to $\M_H(G)$, one has that $H_n(H,
X_f(A/F_{\infty}))$ is finitely generated over $\Op$ for every
$n\geq 0$. This in turn implies that
$\mu_{\Op\ps{\Ga}}\big(H_n(H,X_f(A/F_{\infty}))\big) = 0$ for every
$n\geq 0$, and therefore, it follows that
$\mu_{\Op\ps{G}}\big(X(A/F_{\infty})\big) =
\mu_{\Op\ps{\Ga}}\big(X(A/F_{\cyc})\big)$. \epf

We can now prove Theorem \ref{MHG}.

\bpf[Proof of Theorem \ref{MHG}]
 Let $L$ be a finite extension of $F$ contained in $F_{\infty}$ such
that $G_L=\Gal(F_{\infty}/L)$ is pro-$p$ and has no $p$-torsion.
Since $H_L=\Gal(F_{\infty}/L_{\cyc})$ is a subgroup of $H$ of finite
index, a $\Op\ps{G}$-module is finitely generated over $\Op\ps{H}$
if and only if it is finitely generated $\Op\ps{H_L}$-module.
Therefore, we have that $X(A/F_{\infty})$ belongs to
$\M_{H_L}(G_L)$. The hypothesis of the theorem and Lemma
\ref{torsion implies surjective} allows us to apply Lemma \ref{MHG
lemma} to conclude that $\mu_{\Op\ps{G_L}}\big(X(A/F_{\infty})\big)
= \mu_{\Op\ps{\Ga_L}}\big(X(A/F_{\cyc})\big)$ and
$\mu_{\Op\ps{G_L}}\big(X(A^*/F_{\infty})\big) =
\mu_{\Op\ps{\Ga_L}}\big(X(A^*/F_{\cyc})\big)$. The required
conclusion is now an immediate consequence of Theorem \ref{Zp^r}.
\epf

We end the section stating some further auxiliary results on the
structure of the Selmer group. First, we recall that it is a
well-known observation from the structure theory of a finitely
generated $\Op\ps{\Ga}$-module that a finitely generated
$\Op\ps{\Ga}$-module $M$ is finitely generated over $\Op$ if and
only if $M$ is a torsion $\Op\ps{\Ga}$-module with
$\mu_{\Op\ps{\Ga}}(M)=0$. Such a statement is false when one
replaces $\Ga$ by a general compact $p$-adic Lie group $G$ of
dimension $>1$. For instance, if $G$ has a quotient $\Ga\cong \Zp$,
the module $\Op\ps{\Ga}/p$ is clearly a torsion $\Op\ps{G}$-module
with $\mu_{\Op\ps{G}}(M)=0$, but it is not finitely generated over
$\Op\ps{H}$, where $H$ is the subgroup of $G$ such that $G/H=\Ga$.
In fact, the module $\Op\ps{\Ga}/p$ also belongs to $\M_H(G)$.
Therefore, this example also shows that a $\Op\ps{G}$-module
belonging to $\M_H(G)$ with $\mu_{\Op\ps{G}}(M)=0$ needs not be
finitely generated over $\Op\ps{H}$. However, for the dual Selmer
group, we have the following interesting observation.

\bp
 Suppose that $G$ is a compact pro-$p$ $p$-adic
Lie group with no $p$-torsion. Assume that $X(A/F_{\infty})$ belongs
to $\M_H(G)$. Suppose also that $H^2(G_S(F_{\cyc}), A) = 0$,
$H^2(G_S(F_{\infty}), A) = 0$, and $\la_{A/F_{\cyc}}$ and
$\la_{A/F_{\infty}}$ are surjective. Then $X(A/F_{\infty})$ is
finitely generated over $\Op\ps{H}$ if and only if it belongs to
$\M_H(G)$ with $\mu_{\Op\ps{G}}(X(A/F_{\infty})=0$. \ep

\bpf
 As the ``only if" direction is well-known, we will only prove the ``if"
direction. Suppose that $X(A/F_{\infty})$ belongs to $\M_H(G)$ and
$\mu_{\Op\ps{G}}(X(A/F_{\infty}))=0$. It then follows from Lemma
\ref{MHG lemma} that \linebreak
$\mu_{\Op\ps{\Ga}}(X(A/F_{\cyc}))=0$. By \cite[Proposition
2.5]{CS12}, we have that $X(A/F_{\cyc})$ is $\Op\ps{\Ga}$-torsion.
Thus, one can apply the structure theory of $\Op\ps{\Ga}$-module to
see that $X(A/F_{\cyc})$ is finitely generated over $\Op$. Now
appealing to the argument in Proposition \ref{fg H}, we have that
$X(A/F_{\infty})$ is finitely generated over $\Op\ps{H}$. The proof
of the proposition is now completed. \epf

The next result concerns with comparing the structural properties of
the dual Selmer groups $X(A/F_{\infty})$ and $X(A^*/F_{\infty})$. We
have seen in Proposition \ref{almost abelian} that if $F_{\infty}$
is an almost abelian admissible extension, then $X(A/F_{\infty})$ is
$\Op\ps{G}$-torsion if and only if $X(A^*/F_{\infty})$ is
$\Op\ps{G}$-torsion. Naturally, one may ask if one can establish an
analogous result for the property of belonging to $\M_H(G)$. At this
point of writing, we are only able to establish such result in the
following special case.

\bp Assume that $G=\Zp^2$.  Furthermore, suppose that $A(L_{\cyc})$
and $A^*(L_{\cyc})$ are finite for
 every finite extension $L$ of $F$ contained in
$F_{\infty}$.

Then $X(A/F_{\infty})$ belongs to $\M_H(G)$ if and only if
$X(A^*/F_{\infty})$ belongs to $\M_H(G)$. \ep

\bpf
 It suffices to show that if $X(A/F_{\infty})$ belongs to
 $\M_H(G)$, then $X(A^*/F_{\infty})$ also belongs to $\M_H(G)$.
 Fix a lifting of $\Ga$ to $G$ so that $G = H\times\Ga$.
Let $L_{\infty}$ be the fixed field of the subgroup $\Ga$. Thus,
$L_{\infty}$ is a $\Zp$-extension of $F$ and set $L_n$ to be its
unique subextension of degree $p^n$ over $F$. Now suppose that
$X(A/F_{\infty})$ belongs to $\M_H(G)$. By \cite[Proposition
2.5]{CS12}, we have that $X(A/L_{n, \cyc})$ is $\Op\ps{\Ga}$-torsion
for every $n$. Then the proof of \cite[Theorem 3.8]{CS12} carries
over yielding
\[\mu_{\Op\ps{\Ga}}(X(A/L_{n,\cyc}))
= p^n\mu_{\Op\ps{\Ga}}(X(A/F_{\cyc})).\] By Theorem \ref{Zp^r}, this
in turns implies that $X(A^*/L_{n, \cyc})$ is $\Op\ps{\Ga}$-torsion
for every $n$ and
\[\mu_{\Op\ps{\Ga}}(X(A^*/L_{n,\cyc})) =
p^n\mu_{\Op\ps{\Ga}}(X(A^*/F_{\cyc})).\]
 By the reverse direction of
\cite[Theorem 3.8]{CS12}, this in turn implies that
$X(A^*/F_{\infty})$ belongs to $\M_H(G)$. \epf

\section{Artin twists of Selmer groups} \label{Artin twist of Selmer
groups}
 We retain the notation of the previous section. Let $\rho:
G\lra GL_d(\Op')$ be an Artin representation, where $\Op'$ is the
ring of integers of a finite extension $K'$ of $K$. Denote
$W_{\rho}$ to be a free $\Op'$-module of rank $d$ realizing $\rho$.
The (Greenberg strict) Selmer group attached to the Artin twist of
the data $\big(A, \{A_v\}_{v|p}, \{A^+_v\}_{v|\R} \big)$ is defined
by replacing $A$ by $A\ot_{\Op}W_{\rho}$ and $A_v$ by
$A_v\ot_{\Op}W_{\rho}$ in the definition of the Selmer groups. We
denote this twisted Selmer group by $S(A\ot_{\Op}W_{\rho}/L)$ and
its Pontryagin dual by  $X(A\ot_{\Op}W_{\rho}/L)$.

We should mention that the twisted data $\big(A\ot_{\Op}W_{\rho},
\{A_v\ot_{\Op}W_{\rho}\}_{v|p}, \{(A\ot_{\Op}W_{\rho})^+_v\}_{v|\R}
\big)$ need not satisfy equality (\ref{data equality}). Therefore,
the next two theorems do not follow immediately from the results
before.

\bt \label{Zpr twist}
 Suppose that $F_{\infty}$ is an almost abelian
admissible $p$-adic Lie extension of $F$, and suppose that
$\rho:G\lra GL_d(\Op')$ is an Artin representation with
$\Op\subseteq \Op'$.
 Then the dual Selmer groups $X(A\ot_{\Op}W_{\rho}/F_{\infty})$ and
 $X(A^*\ot_{\Op}W_{\hat{\rho}}/F_{\infty})$ have
the same
 $\Op'\ps{G}$-rank and the same $\mu_{\Op'\ps{G}}$-invariant.
\et

\bt \label{MHG twist} Let $F_{\infty}$ be an $S$-admissible $p$-adic
Lie extension. Assume that both $X(A/F_{\infty})$ and
$X(A^*/F_{\infty})$ belong to $\M_H(G)$. Furthermore, suppose that
$A(L_{\cyc})$ and $A^*(L_{\cyc})$ are finite for
 every finite extension $L$ of $F$ contained in
$F_{\infty}$.

For each Artin representation $\rho:G\lra GL_d(\Op')$ with
$\Op\subseteq \Op'$, we have
\[ \mu_{\Op'\ps{G}}\big( X(A\ot_{\Op}W_{\rho}/F_{\infty})\big) =
\mu_{\Op'\ps{G}} \big(X(A\ot_{\Op}W_{\hat{\rho}}/F_{\infty})\big).
\] \et

We record the following lemma.

\bl \label{twist lemma}
 For each Artin representation $\rho:G\lra GL_d(\Op')$, we have
  \[ \mathrm{tw}_{\hat{\rho}}\big(X(A/F_{\infty})\big)
  = X(A\ot_{\Op}W_{\rho}/F_{\infty}),  \]
  where $\hat{\rho}$ is the contragredient representation of $\rho$.
\el

\bpf
 For any $\Op$-module $M$, we write $M_{\Op'} = M \ot_{\Op}
\Op'$. Since $\rho$ factors through a finite quotient of $G$, we
have $S(A\ot_{\Op}W_{\rho}/F_{\infty}) = S(A/F_{\infty})_{\Op'}
\ot_{\Op'}W_{\rho}$. Hence
\[ X(A\ot_{\Op}W_{\rho}) = \Hom_{\Op'}(W_{\rho},
X(A/F_{\infty})_{\Op'}) = \Hom_{\Op'}\big(\Hom_{\Op'}(\Op',
W_{\hat{\rho}}),X(A/F_{\infty})_{\Op'}\big).
  \]
 Since $W_{\hat{\rho}}$ is a free $\Op'$-module, it follows from
\cite[Chapter VI, Proposition 5.2]{CE} that the latter module is
isomorphic to
 \[\Hom_{\Op'}(\Op',
X(A/F_{\infty})_{\Op'})\ot_{\Op'}W_{\hat{\rho}} =
\mathrm{tw}_{\hat{\rho}}\big(X(A/F_{\infty})\big). \] This proves
the lemma. \epf

The conclusion of Theorem \ref{Zpr twist} is now plain from Lemmas
\ref{mu lemma Artin twist} and \ref{twist lemma}, and Theorem
\ref{almost abelian}. The conclusion of Theorem \ref{MHG twist} will
follow from Lemmas \ref{mu lemma Artin twist} and \ref{twist lemma},
and Theorem \ref{MHG}.

\section{Comparing Greenberg Selmer groups and Selmer complexes}
\label{Greenberg Selmer
groups}

As before, let $p$ be a fixed prime. Furthermore, we assume that $F$
has no real primes if $p=2$. Let $\big(A, \{A_v\}_{v|p},
\{A^+_v\}_{v|\R}\big)$ denote the data defined in Section
\ref{compare} and $\big(A^*, \{A^*_v\}_{v|p},
\{(A^*)^+_v\}_{v|\R}\big)$ the dual data.  Fix a finite set $S$ of
primes of $F$ which contains all the primes above $p$, the ramified
primes of $A$ and the infinite primes. Now set
\[ H^1_{Gr}(F_v, A)=
\begin{cases} \ker\big(H^1(F_v, A)\lra H^1(F_v, A/A_v)\big) & \text{\mbox{if} $v|p$},\\
 \ker\big(H^1(F_v, A)\lra H^1(F^{ur}_v, A)\big) & \text{\mbox{if} $v\nmid p$.}
\end{cases} \]

The Greenberg Selmer group attached to these data is then defined by
\[ \Sel^{Gr}(A/F) = \ker\Big( H^1(G_S(F),A)\lra \bigoplus_{v \in S}H^1_g(F_v,
A)\Big),\] where we write $H^1_g(F_v, A) = H^1(F_v, A)/H^1_{Gr}(F_v,
A)$. For an $S$-admissible $p$-adic Lie extension $F_{\infty}$, we
define $\Sel^{Gr}(A/F_{\infty}) = \ilim_L \Sel^{Gr}(A/L)$ and denote
$X^{Gr}(A/F_{\infty})$ to be the Pontryagin dual of
$\Sel^{Gr}(A/F_{\infty})$. The following lemma compares the two
Greenberg Selmer groups.

\bl
  We have an exact sequence
\[0\lra S(A/F_{\infty})\lra \Sel^{Gr}(A/F_{\infty}) \lra N\lra 0,\]
where $N$ is a cofinitely generated $\Op\ps{H}$-module. \el

\bpf
 Now consider the following commutative
diagram
\[  \entrymodifiers={!! <0pt, .8ex>+} \SelectTips{eu}{}\xymatrix{
    0 \ar[r]^{} & S(A/F_{\infty}) \ar[d] \ar[r] &  H^1(G_S(F_{\infty}), A) \ar@{=}[d]
    \ar[r] & \bigoplus_{v \in S} J_v(A/F_{\infty}) \ar[d]^{\al} \\
    0 \ar[r]^{} & \Sel^{Gr}(A/F_{\infty}) \ar[r]^{} & H^1(G_S(F_{\infty}), A) \ar[r] & \
    \bigoplus_{v \in S} J^{Gr}_v(A/F_{\infty})  } \]
with exact rows, where $J^{Gr}_v(A/F_{\infty}) =
\ilim_L\bigoplus_{w|v} H^1_g(L_w, A)$. It therefore remains to show
that $\ker \al$ is  cofinitely generated over $\Op\ps{H}$. Clearly,
$J_v(A/F_{\infty}) = J^{Gr}_v(A/F_{\infty})$ for $v\nmid p$. For
$v|p$, choose a prime $w$ of $F_{\infty}$ above $v$. Write
$I_{\infty, w}$ for the inertia subgroup of
$\Gal(\overline{F}_{\infty,w}/F_{\infty, w})$ and $U_w =
\Gal(\overline{F}_{\infty,w}/F_{\infty, w})/I_{\infty, w}$. It then
follows from the Hochschild-Serre spectral sequence that we have
\[ 0\lra H^1(U_w, (A/A_v)^{I_{\infty,w}}) \lra H^1(H_w, A/A_v)
\lra H^1(I_{\infty,w}, A/A_v)^{U_w}.\]
 Since $U_w$ is topologically cyclic, $H^1(U_w,
(A/A_v)^{I_{\infty,w}}) \cong
\big((A/A_v)^{I_{\infty,w}}\big)_{U_v}$ and so is cofinitely
generated over $\Op$. Now since the decomposition group of $G$ at
$v$ has at least dimension one for each $v|p$, it follows that $\ker
\al$ is cofinitely generated over $\Op\ps{H}$, as required. \epf

\bl \label{Greenberg equal} One has  \[ \rank_{\Op\ps{G}}\big(
X(A/F_{\infty})\big) = \rank_{\Op\ps{G}}
\big(X^{Gr}(A/F_{\infty})\big)
\] and \[ \mu_{\Op\ps{G}}\big( X(A/F_{\infty})\big) =
\mu_{\Op\ps{G}} \big(X^{Gr}(A/F_{\infty})\big). \]
 Furthermore, $X(A/F_{\infty})$ belongs to $\M_H(G)$ if and only
if $X^{Gr}(A/F_{\infty})$ belongs to $\M_H(G)$.  \el

\bpf
 By the preceding lemma, one has an exact sequence
 \[0\lra N'\lra X^{Gr}(A/F_{\infty}) \lra X(A/F_{\infty})\lra 0\]
for some finitely generated $\Op\ps{H}$-module $N'$. The first
equality and the final assertion are then immediate, and the second
equality follows from \cite[Lemma 2.1(c)(iii)]{Lim2}. \epf

\br
 Note that Lemma \ref{Greenberg equal} does not require any
 torsion assumptions on either of the Selmer groups.
\er

The next two theorems are immediate from combining the results in
Sections \ref{compare} and \ref{compare2} with Lemma  \ref{Greenberg
equal}.

\bt
 Suppose that $F_{\infty}$ is an almost abelian
admissible $p$-adic Lie extension of $F$.
 Then $X^{Gr}(A/F_{\infty})$ and $X^{Gr}(A^*/F_{\infty})$ have
the same
 $\Op\ps{G}$-rank and the same $\mu_{\Op\ps{G}}$-invariant.
\et

\bt  Let $F_{\infty}$ be an $S$-admissible $p$-adic Lie extension.
Assume that both $X^{Gr}(A/F_{\infty})$ and $X^{Gr}(A^*/F_{\infty})$
belong to $\M_H(G)$.  Furthermore, suppose that $A(L_{\cyc})$ and
$A^*(L_{\cyc})$ are finite for
 every finite extension $L$ of $F$ contained in
$F_{\infty}$.

 Then we have
\[ \mu_{\Op\ps{G}}\big( X^{Gr}(A/F_{\infty})\big) =
\mu_{\Op\ps{G}} \big(X^{Gr}(A^*/F_{\infty})\big). \] \et

\bigskip We now consider the Selmer complex of $\big(A, \{A_v\}_{v|p},
\{A^+_v\}_{v|\R}\big)$. The notion of a Selmer complex was first
conceived and introduced in \cite{Nek}. In our discussion, we will
considered a modified version of Selmer complexes as given in
\cite[4.2.11]{FK}. Write $T^* = \Hom_{\cts}(A,\Qp/\Zp)$ and $T^*_v =
\Hom_{\cts}(A/A_v,\Qp/\Zp)$. For every finite extension $L$ of $F$
and $w$ a prime of $L$ above $p$, write $T^*_w = T^*_v$, where $v$
is the prime of $F$ below $w$. For any profinite group $\mathcal{G}$
and a topological abelian group $M$ with a continuous
$\mathcal{G}$-action, we denote by $C(\mathcal{G}, M)$ the complex
of continuous cochains of $\mathcal{G}$ with values in $M$. Let
$F_{\infty}$ be an $S$-admissible extension of $F$ with Galois group
$G$. We define a $(\Op\ps{G})[G_S(F)]$-module $\mathcal{F}_G(T^*)$
as follows: as a $\Op$-module, $\mathcal{F}_G(T^*) =
\Op\ps{G}\ot_{\Op}T^*$, and the action of $G_S(F)$ is given by the
formula $\sigma(x\ot t) = x\bar{\sigma}^{-1}\ot \sigma t$, where
$\bar{\sigma}$ is the canonical image of $\sigma$ in $G \subseteq
\Op\ps{G}$. We define the $(\Op\ps{G})[\Gal(\bar{F}_v/F_v)]$-module
$\mathcal{F}_G(T_v^*)$ in a similar fashion.

For every prime $v$ of $F$, we write $C\big(F_v,
\mathcal{F}_G(T^*)\big)= C\big(\Gal(\bar{F}_v/F_v),
\mathcal{F}_G(T^*)\big)$. For each prime $v$ not dividing $p$, let
$C_f\big(F_v, \mathcal{F}_G(T^*)\big)$ be the subcomplex of
$C\big(F_v, \mathcal{F}_G(T^*)\big)$, whose degree $m$-component is
$0$ unless $m\neq  0, 1$, whose degree $0$-component is
$C^0\big(F_v, \mathcal{F}_G(T^*)\big)$, and whose degree
$1$-component is
\[ \ker\Big( C^1\big(F_v, \mathcal{F}_G(T^*)\big)_{d=0} \longrightarrow
H^1\big(F_v^{ur},\mathcal{F}_G(T^*)\big) \Big). \]

The Selmer complex $SC(T^*, T^*_v)$ is defined to be
\[ \mathrm{Cone}\bigg( C\big(G_S(F), \mathcal{F}_G(T^*)\big) \longrightarrow
\displaystyle\bigoplus_{v|p}C\big(F_v,
\mathcal{F}_G(T^*)/\mathcal{F}_G(T^*_v)\big)\oplus \displaystyle
\bigoplus_{v\nmid p} C\big(F_v, \mathcal{F}_G(T^*)\big)/C_f\big(F_v,
\mathcal{F}_G(T^*)\big)\bigg)[-1].
\]
 Here $[-1]$ is the translation by $-1$ of the complex.
We now state the following proposition which is given in
\cite[Proposition 4.2.35]{FK}.

\bp Let $\mathcal{G}$ be the kernel of $Gal(\bar{F}/F) \lra G$. For
a place $v$ of $F$, fixing an embedding $F \hookrightarrow F_v$, let
$\mathcal{G}(v)$ be the kernel of $Gal(\bar{F}_v/F_v)\lra G$ and let
$G_v\subseteq G$ be the image. Then the following statements hold.
\begin{enumerate}
\item[$(a)$] $H^i\big(SC(T^*, T^*_v)\big) = 0$ for $i\neq 1, 2,
3$.

\item[$(b)$] We have an exact sequence
\[ \ba{c}
0 \lra X(A/F_{\infty}) \lra H^2\big(SC(T^*, T^*_v)\big) \lra
\displaystyle \bigoplus_{v|p} \Op\ps{G}\ot_{\Op\ps{G_v}}
\big(T^*_v(-1)\big)_{\mathcal{G}(v)}
 \hspace{2in}\\
 \hspace{1.5 in} \lra \big(T^*(-1)\big)_{\mathcal{G}} \lra H^3\big(SC(T^*,
 T^*_v)\big) \lra 0.
 \ea\]
\end{enumerate}
\ep

\medskip
Denote $SC(T, T_v)$ to be the Selmer complex of the dual data
$\big(A^*, \{A^*_v\}_{v|p},  \{(A^*)^+_v\}_{v|\R}\big)$. We can now
state the following analogous results for the second cohomology
groups of the Selmer complexes.

\bt
 Suppose that $F_{\infty}$ is an almost abelian
admissible $p$-adic Lie extension of $F$.
 Then $ H^2\big(SC(T, T_v)\big)$ and $ H^2\big(SC(T^*, T^*_v)\big)$ have
the same
 $\Op\ps{G}$-rank and the same $\mu_{\Op\ps{G}}$-invariant.
\et

\bt  Let $F_{\infty}$ be an $S$-admissible $p$-adic Lie extension of
dimension $>1$. Assume that both $X(A/F_{\infty})$ and
$X(A^*/F_{\infty})$ belong to $\M_H(G)$.  Furthermore, suppose that
$A(L_{\cyc})$ and $A^*(L_{\cyc})$ are finite for
 every finite extension $L$ of $F$ contained in
$F_{\infty}$.

 Then we have
\[ \mu_{\Op\ps{G}}\big(  H^2\big(SC(T, T_v)\big)\big) =
\mu_{\Op\ps{G}} \big(  H^2\big(SC(T^*, T^*_v)\big)\big). \] \et

\section{Examples} \label{examples}

We discuss some arithmetic examples of our main results. For
brevity, we will only consider the strict Selmer groups. Our choices
of examples are motivated by the possible application of Theorem
\ref{MHG}. Of course, our Theorem \ref{Zp^r} is unconditional and
can be applied to many other examples that are not discussed here.
We refer readers to \cite[Section 1.2]{We} for some of these
examples.

\subsection{Abelian varieties}

Let $B$ be an abelian variety of dimension $g$ defined over a number
field $F$. For simplicity, we will assume that $F$ is totally
imaginary and that $B$ has semistable reduction over $F$. We define
a data $(A, \{A_v\})$ by first setting $A = B[p^{\infty}]$. For each
prime $v$ of $F$ above $p$, let $\mathcal{F}_v$ be the formal group
attached to the Neron model for $B$ over the ring of integers
$\Op_{F_v}$ of $F_v$, and we assume that $\mathcal{F}_v$ is a formal
group of height $g$ for all $v|p$. For instance, this is satisfied
if $B$ has good ordinary reduction at all $v|p$. We then set $A_v =
\mathcal{F}_v(\overline{\m})[p^{\infty}]$, where $\overline{\m}$ is
the maximal ideal of the rings of integers of $\overline{F}_v$. Note
that $A_v\cong (\Qp/\Zp)^g$ as an abelian group by our height
assumption.

Then the dual data $(A^*, \{A^*_v\})$ is given by $A^* =
B^t_{p^{\infty}}$ and $A^*_v =
\mathcal{F}^t_v(\overline{\m})[p^{\infty}]$, where $B^t$ is the dual
abelian variety of $B$ and $\mathcal{F}^t_v$ is the formal group
attached to the Neron model for $B^t$ over the ring of integers
$\Op_{F_v}$ of $F_v$. Clearly, $\mathcal{F}^t_v$ has height $g$ and
$A^*_v\cong (\Qp/\Zp)^g$ as an abelian group. It is worthwhile to
mention that the Greenberg strict Selmer groups attached to these
data coincide with the classical Selmer groups of the abelian
variety (see \cite{CG}), when the Selmer groups are considered over
an admissible $p$-adic Lie extension.

Now if $F_{\infty}$ is an almost abelian $S$-admissible extension of
$F$, it follows from Theorem \ref{almost abelian} that
$X(A/F_{\infty})$ and $X(A^*/F_{\infty})$ have the same
$\Zp\ps{G}$-rank and the same $\mu_{\Zp\ps{G}}$-invariant. In
particular, this applies when $B$ is an abelian variety with complex
multiplication and $F_{\infty} = F(B[p^{\infty}])$.

We now consider the situation of a general noncommutative $p$-adic
Lie extension. The conditions on the finiteness of $A(L^{\cyc})$ and
$A^*(L^{\cyc})$ are consequences of \cite[Theorem 4.3]{Win}.
Therefore, to apply Theorem \ref{MHG}, it remains to require that
$X(A/F_{\infty})$ and $X(A^*/F_{\infty})$ both belong to $\M_H(G)$.
When $B$ is an elliptic curve, this is conjectured in \cite{CFKSV}.
In view of their conjecture, it seems quite reasonable to expect
that $X(A/F_{\infty})$ and $X(A^*/F_{\infty})$ both belong to
$\M_H(G)$ for a general abelian variety $B$ which has good ordinary
reduction at all primes above $p$. However, we should mention that
there is currently no known method to prove that the dual Selmer
group belongs to $\M_H(G)$ in general.

 We like to mention that in a preprint of Bhave \cite{Bh}, she
was able to establish the equality of $\mu$-invariant for the case
that $B$ is an abelian variety without complex multiplication under
a weaker torsion condition on the Selmer groups than our Theorem
\ref{MHG}. Therefore, our results in the abelian varieties cases
complement the results there.

We can also apply our results (Theorems \ref{Zpr twist} and \ref{MHG
twist}) to the Artin twists of the abelian variety $B$.

\subsection{Hilbert modular forms}
Let $F$ be totally real number field and let $f$ be a primitive
Hilbert modular form of parallel weight $(k, k, ..., k)$ with $k
\geq 2$ which is $p$-ordinary. Denote $K_f$ to be the field
generated by the coefficient of $f$. We will write $K$ for the
localization of $K_f$ at some fixed prime of $K_f$ above $p$. Let
$V$ be the two-dimensional representation over $K$ associated to
$f$. Since $f$ is assumed to be $p$-ordinary, it follows that for
every prime $v$ of $F$ above $p$, we have a one-dimensional subspace
$V_v$ of $V$ which is $\Gal(\bar{F}_v/F_v)$-invariant and has the
property that the inertia subgroup acts by a power of the cyclotomic
character on this subspace and trivially on the one dimensional
quotient. Fix a Galois lattice $T$ of $V$, and for each prime $v$ of
$F$ above $p$, we set $T_v = T\cap V_v$. Then in this situation, the
data $(A, \{A_v\}_{v|p})$ is given by $A= V/T$ and $A_v = V_v/T_v$.
The dual data is then the above definition for the dual modular form
$\bar{f}$. In view of Lemma \ref{data base change}, we can apply our
results to an admissible $p$-adic Lie extension $F_{\infty}$ of $F$
which is either totally real or totally imaginary. Again, if
$F_{\infty}$ is an almost abelian $S$-admissible extension of $F$,
Theorem \ref{almost abelian} applies.

We now consider the case of a general noncommutative $p$-adic Lie
extension and suppose that $f$ is a primitive cuspidal modular form
of positive weight $k\geq 2$ with $K_f=\Q$. The finiteness of
$A(L^{\cyc})$ and $A^*(L^{\cyc})$ are established in the proof of
\cite[Lemma 2.2]{Su}. In order to apply Theorem \ref{MHG}, one will
require that $X(A/F_{\infty})$ and $X(A^*/F_{\infty})$ both belong
to $\M_H(G)$. This latter condition will follow from the $\M_H(G)$
conjecture for Hida families (see \cite{CS12}). However, it does not
seem to be known how one can verify the $\M_H(G)$-condition in
either situations.

Finally, we mention that we can apply our results to Artin twists of
the two-dimensional representation of $f$ in a similar fashion as in
the previous subsection.

\end{document}